\documentclass[11 pt]{article}
\usepackage{amssymb,amsmath}
\usepackage{amsfonts}
\usepackage{a4wide}
\usepackage[latin1]{inputenc}

\newtheorem{coro}{Corollary}
\newtheorem{nota}{Notations}
\newtheorem{defi}{Definition}
\newtheorem{lem}{Lemma}
\newtheorem{claim}{Claim}
\newtheorem{prop}{Proposition}
\newtheorem{theo}{Theorem}
\newtheorem{rem}{Remark}

\newenvironment{demo}
{\medbreak\noindent{\sc Proof :}}
{\hfill$\square$\medbreak}
\newenvironment{demof}[1]
{\medbreak\noindent{\sc Proof of {#1} :}}
{\hfill$\square$\medbreak}

\begin{document}

\newcommand\C{{\mathcal C}}

\title{A Proof of Gromov's Algebraic Lemma}
\author{David Burguet}

\maketitle

\pagestyle{myheadings} \markboth{\normalsize\sc David
Burguet}{\normalsize\sc A Proof of Gromov's Algebraic Lemma}

\noindent
\textbf{Abstract :} Following the analysis of differentiable
mappings of Y. Yomdin, M. Gromov  has stated a very elegant
``Algebraic Lemma'' which says that the ``differentiable size'' of
an algebraic subset may be bounded in terms only of its dimension,
degree and diameter - regardless of the size of the underlying coefficients.
 We give a complete and elementary proof of Gromov's result using
the ideas presented in his Bourbaki talk as well as other
necessary ingredients.

\section{Introduction}
Several problems in, e.g., Analysis and Dynamical Systems, require
estimating the differentiable size of semi-algebraic subsets. Y. Yomdin
developped many tools to this end \cite{C&Y}. M. Gromov observed that one of these tools could be refined to give the following very elegant statement :

\begin{theo} For all integers $r\geq 1$, $d\geq 0$, $\delta\geq 0$, 
there exists $M<\infty$ with the following properties.
For any semi-algebraic compact subset $A\subset ]0,1[^d$ of maximum
dimension $l$ and of degree $\leq \delta$, there exist an integer
$N$ and maps
$\phi_1,...,\phi_N:[0,1]^l \mapsto ]0,1[^d$ satisfying 
 $\bigcup _{i=0}^N \phi_i([0,1]^l)=A$, such that : 

\begin{itemize}
\item $\|\phi_{i/]0,1[^l}\|_r:=\max _{\beta:|\beta|\leq r} \|
\partial ^{\beta} \phi_i\|_{\infty}\leq1$ ;
\item $N\leq M$ ;
\item $deg(\phi_i)\leq M$.
\end {itemize}

\end{theo}

In his S\'eminaire Bourbaki \cite{Gr}, M. Gromov gives many ideas but stops short of a complete proof. On the other hand, this result
has been put to much use, especially in Dynamical System Theory. Y. Yomdin
 \cite{Yoma},\cite{Yomb} used it to compare the topological entropy
and the ``homological size''  for $\C^r$ maps (in
particular, Y. Yomdin
 proves in \cite{Yoma}  Shub's conjecture in the case of
$\C^{\infty}$ maps). S. Newhouse \cite{New} then showed, using Pesin's
theory, how this gives, for $\C^{\infty}$ smooth maps, upper-semicontinuity
of the metric entropy and therefore the existence of invariant
measures with maximum entropy. J. Buzzi \cite{Bu} observed that in fact
Y. Yomdin's
estimates give a more uniform result called asymptotic h-expansiveness,
which was in turn used by M. Boyle, D. Fiebig and U. Fiebig \cite{BFF} to
prove existence of principal symbolic extensions. The dynamical consequences
of the above theorem are still developping in the works of M. Boyle, T. Downarowicz, S. Newhouse  and others \cite{New1},\cite{BD}.\\

The proof of this theorem is trivial in dimension 1 and easy
in dimension 2 (see part 6). To prove the theorem
in higher dimensions, we introduce the notion of triangular
$(\C^{\alpha},K)$-Nash maps : it is the subject of the part 3. 
Part 4 is devoted to the structure of semi-algebraic sets.
In  part 5, by taking the limit of ``good'' parametrizations,
 we reduce the main theorem to a proposition about the parametrization of
semi-algebraic "smooth" maps (thus avoiding the singularities). The other
difficulties are dealt with as suggested by M. Gromov.
 The proof by induction of this
proposition is done in the last section. Describe briefly the structure of this
proof. We distinguish three independent steps : 

\begin{itemize}
\item we consider a semi-algebraic map defined on a subset of higher
dimension and we bound the first derivative in the first coordinate.
\item we bound the derivative of higher order in the first coordinate.
\item fixing the dimension of the semi-algebraic set and the order of 
derivation, we bound the next derivative for the order defined on 
$\mathbb{N}^d$ in part 3.  
\end{itemize}

\medbreak
As I was completing the submission of this paper, I learnt that A. Wilkie
had written a proof of the same theorem \cite{w}. I am grateful to M. Coste
for this reference. In the first version of this article, M. Coste also pointed out  a mistake corrected here by Remark \ref{adapte}. 
\medbreak

\section{Semi-algebraic sets and maps}

First recall some basic results concerning semi-algebraic sets.
We borrow from \cite{Coste}. For completeness, other references are
\cite{BR},\cite{BCR},\cite{C&Y}.

\begin{defi} $A\subset \mathbb{R}^d$ is a \underline{semi-algebraic set} if it
can be written as a finite union of sets of the form
 $\{x\in \mathbb{R}^d \ | \ P_1(x)>
0,...P_r(x)> 0, P_{r+1}(x)=0,...,P_{r+s}(x)=0\}$, where $r,s\in
\mathbb{N}$ and $P_1,...,P_{r+s}\in \mathbb{R}[X_1...X_d]$ . Such a
formula is called a \underline{presentation} of $A$.

The degree of a presentation is the sum of the total degrees of the
polynomials involved (with multiplicities). The \underline{degree} of a
semi-algebraic set is the minimum  degree of its presentations.
\end{defi}

\begin{defi}\label{sam} $f:A\subset \mathbb{R}^d\rightarrow
\mathbb{R}^n$ is a  \underline{semi-algebraic map} if the graph of $f$ is a
semi-algebraic set.
\end{defi}

\begin{defi}\label{nashm} A \underline{Nash manifold} is an analytic submanifold of $\mathbb{R}^d$,
which is a semi-algebraic set.

 A \underline{Nash map} is a map defined on a Nash manifold, which is analytic
and semi-algebraic.
\end{defi}

We have the following description of a semi-algebraic set (See
\cite{Coste}, Prop. 3.5 p 124 and see \cite{C&Y} Prop. 4.4 p 48) :

\begin{theo}\label{thm}(stratification) Let $A\subset
\mathbb{R}^n$ be a semi-algebraic set. There exist an
integer $N$ (bounded in terms of $deg(A)$) and connected Nash
manifolds $A_1,...,A_N$ such that $A=\coprod_{i=1}^N A_i$ and
$\forall j\neq i$ $(A_i\bigcap adh(A_j)\neq\emptyset)\Rightarrow
(A_i\subset adh(A_j) \
 et \ dim(A_i)<dim(A_j))$. ($\coprod$ : disjoint union).
\end{theo}

\begin{defi} In the notations of the previous proposition, the 
\underline{maximum dimension} of $A$  is the maximum dimension of the Nash
manifolds $A_1,...A_N$.
\end{defi}

\section{ $(\C^{\alpha},K)$-Nash maps and triangular maps}\label{app}

\begin{defi}\label{def}  $\mathbb{N}^d$ is provided with the
order
$\preceq$, defined as follows : 

for $\alpha=(\alpha_1,...\alpha_d)$, $\beta=(\beta_1,...\beta_d)\in
\mathbb{N}^d$

$\alpha\preceq \beta$ iff ($|\alpha|:=\sum_i \alpha_i<|\beta|$) or
($|\alpha|=|\beta|$ et $\alpha_k\leq\beta_k$, where $k:=\max\{l\leq
n \ :\ \alpha_l\neq\beta_l\}$)
\end{defi}

\begin{nota} The order $\preceq$ is a total order. Hence, for $\alpha \in \mathbb{N}^d$, we can set : $$\alpha +1:=\min \{\beta \in \mathbb{N}^d : \alpha \preceq \beta \ \ and \  \alpha \neq \beta\}$$.
\end{nota}

\begin{defi}\label{def} Let $K\in\mathbb{R}^{+}$, $d\in \mathbb{N}$,
$\alpha\in\mathbb{N}^{d}-\{0\}$. Let $A\subset ]0,1[^d$ be a
semi-algebraic open set. 

A map $f:A\rightarrow
]0,1[^d$ is a \underline{$\C^0$-Nash map}, if $f:=(f_1,...f_d)$ is a Nash map,
which can be continuously extended to $adh(A)$. We call again $f$
this unique extension. 

A map $f:A\rightarrow ]0,1[^d$ is a
\underline{$(\C^{\alpha},K)$-Nash map}, if $f$ is a  $\C^0$-Nash map and if
$\|f\|_{\alpha}:=max _{\beta \preceq \alpha,1\leq i\leq d}
\|\partial ^{\beta} f_i\|_{\infty}\leq K$. If $\alpha=(0,0...,r)$,
we write $(\C^r,K),\|.\|_r$ instead of $(\C^{\alpha},K),\|.\|_{\alpha}$.
\end{defi}

The two following lemmas deal with the composition of
$(\C^{\alpha},1)$-Nash maps.

\begin{lem}\label{comp} For all $d,r\in \mathbb{N}^{*}$, there exists  $K<+\infty$,
such that if $\psi,\phi:]0,1[^d \rightarrow ]0,1[^d$ are two
 $(\C^r,1)$-Nash maps, then $\psi \circ \phi$ is a
$(\C^r,K)$-Nash map.

\end{lem}

\begin{demo} Immediate.
\end{demo}

One of the key points of the proof of Gromov's lemma is to control the
 derivatives one after one. This is made possible by the folllowing definition.

\begin{defi}\label{def} We say that a map
$\psi:]0,1[^l\rightarrow ]0,1[^d$ is \underline{triangular} if $l\leq d$ and if there
exists a family of maps $(\psi_i:]0,1[^{min(l,d+1-i)} \rightarrow
]0,1[)_{i=1...d}$, such that
$$\psi=(\psi_1(x_1...x_l),...,\psi_{d-l+1}(x_1...x_l),\psi_{d-l+2}(x_2...x_l),...
,\psi_{d-l+k}(x_k...x_l),...,\psi_d(x_l))$$
\end{defi}

\begin{rem} If $\psi:]0,1[^n\rightarrow ]0,1[^m$ and
$\phi:]0,1[^m\rightarrow ]0,1[^p$ are  triangular, then so is 
$\phi\circ\psi:]0,1[^n\rightarrow ]0,1[^p$.
\end{rem}

In the case of triangular maps, we give the following version of the
lemma \ref{comp}. This result allows an induction on $\alpha\in\mathbb{N}^d$
rather than $r\in \mathbb{N}$, in the proof of the proposition 4.

\begin{lem}\label{lemcomp} For all $d,r\in \mathbb{N}^{*}$, there exists $K<+\infty$
such that if  $\psi,\phi:]0,1[^d \rightarrow ]0,1[^d$ are two
triangular $(\C^{\alpha},1)$-Nash maps with $|\alpha|=r$, then
$\psi\circ\phi$ is a $(\C^{\alpha},K)$-Nash map.
\end{lem}

\begin{demo} Immediate.
\end{demo}

\begin{defi}\label{def} (\underline{resolution of a semi-algebraic set}) Let $M:\mathbb{N}^3\rightarrow \mathbb{R}^+$
and let $K\in \mathbb{R}^+$, $d\in\mathbb{N}^{*}$. Let $A\subset
[0,1]^d$ be a semi-algebraic set of maximum dimension $l$ and
let $\alpha\in\mathbb{N}^{d}-\{0\}$. The family of maps
$(\phi_i:]0,1[^l\rightarrow]0,1[^{d})_{i=1...N}$  
 is a $(M)$-resolution [resp. $(\C^{\alpha},K,M)$-resolution] of $A$ if :
\begin{itemize}
           \item each $\phi_i$ is triangular ;
           \item each $\phi_i$ is a Nash map \footnote{not necessarily $\C^0$-Nash map}  [resp. a
 $(\C^{\alpha},K)$-Nash map] ;
           \item $A=\bigcup_{i=1}^N \phi_i(]0,1[^l)$ \footnote{by convention $]0,1[^{0}=\{0\}$}  [resp. $adh(A)=\bigcup_{i=1}^N \phi_i([0,1]^l)$]  
           \item $N$, $deg(\phi_i)$ are less than $M(0,d,deg(A))$
           [resp.
           $M(|\alpha|,d,deg(A))$].
         \end{itemize}
\end{defi}

\begin{defi} (\underline{resolution of a family of maps})
 Let $M:\mathbb{N}^4\rightarrow \mathbb{R}^+$, $K\in \mathbb{R}^+$,
 $d\in\mathbb{N}^{*}$, $\alpha\in\mathbb{N}^{d}-\{0\}$ and let $f_1,...,f_k:A\rightarrow ]0,1[$ be semi-algebraic maps, where $A\subset ]0,1[^d$ is a semi-algebraic set of maximum dimension $l$. The family of maps $(\phi_i:]0,1[^l\rightarrow
]0,1[^d)_{i=1,...,N}$  is
 a $(\C^{\alpha},K,M)$-resolution of $(f_i)_{i=1,...,k}$ if :

\begin{itemize}
    \item each $\phi_i$ is triangular ;
    \item each $\phi_i$, $f_j\circ\phi_i$ is a  $(\C^{\alpha},K)$-Nash map ;
    \item $adh(A)=\bigcup_{i=1}^N \phi_i([0,1^l])$ ;
    \item $N$, $deg(\phi_i)$, $deg(f_j\circ\phi_i)$ are less than 
 $M(|\alpha|,d,k,max_j(deg(f_j)))$.
\end{itemize}
\end{defi}

We shall consider only functions $M$ in the above setting that are 
independent of the algebraic datas (\emph{i.e.} the functions $f_1,...,f_k$ or the set $A$). Such function can be called ``universal''. By a 
 $(\C^{\alpha},K)$-resolution, we mean a $(\C^{\alpha},K,M)$-resolution with a universal function $M$.  

The following remark is very useful later on :

\begin{lem}\label{lem3}  For all  $M:\mathbb{N}^2\rightarrow \mathbb{R}^+$, there exists $M':\mathbb{N}^2\rightarrow \mathbb{R}^+$ such that we have the following property.

Let  $d\in \mathbb{N}^{*}$,
 $\alpha\in \mathbb{N}^d-\{0\}$. If $f:]0,1[^d\rightarrow ]0,1[$ is a $(\C^{\alpha},M(|\alpha|,d))$-Nash map,
  then $f$ admits a $(\C^{\alpha},1,M')$-resolution.
\end{lem}

\begin{demo} Linear reparametrizations.
\end{demo}

\section{Tarski's Principle}\label{princ}

\begin{prop} (Tarski's principle) Let $A$ a semi-algebraic set of
$\mathbb{R}^{d+1}$ and $\pi:\mathbb{R}^{d+1}\rightarrow\mathbb{R}^d$
the projection defined by $\pi(x_1,...,x_{d+1})=(x_1,...,x_d)$ then
$\pi(A)$ is a semi-algebraic set and $deg(\pi(A))$ is bounded by a
function of $deg(A)$ and $d$.
\end{prop}

\begin{demo} See \cite{BCR} Thm 2.2.1, p 26 and \cite{C&Y} Prop 4.3 p 48
\end{demo}

\begin{coro}\label{quant} Any formula combining sign conditions on semi-algebraic functions by conjonction, disjunction,
negation and universal and existential quantifiers defines a
semi-algebraic set.
\end{coro}

\begin{coro}\label{cor} Let $f:A\subset \mathbb{R}^k \rightarrow \mathbb{R}^l$ be a semi-algebraic function,
then
 $A$ is a
semi-algebraic set and $deg(A)$ is bounded by a function of
$deg(f),k$ and $l$.
\end{coro}

\begin{coro}\label{cor3} If $\phi$ and $\psi$ are two semi-algebraic maps, 
such that the composition $\phi\circ\psi$ is well defined,
 then  $\phi\circ\psi$ is a semi-algebraic map and its degree is bounded by 
a function of $deg(\phi)$ and $deg(\psi)$.
\end{coro}

\begin{demo} See \cite{BCR} Prop 2.2.6 p 28
\end{demo}

\begin{prop}\label{pro1}
For all  $P_1,...,P_s\in \mathbb{R}[X_1,...,X_{d+1}]$, there exist a
partition of $]0,1[^d$ into Nash manifolds $\{A_1,...,A_m\}$ and a finite
family of  Nash maps,
$\zeta_{i,1}<...<\zeta_{i,q_i}$
$:A_i\rightarrow ]0,1[$, for all $1 \leq i \leq m$, such that :
\begin{itemize}
    \item for each $i$ and each $k$, the sign  $P_k(x_1,y)$,
    with $x_1\in ]0,1[$ et $y:=(x_2,...,x_{d+1})\in A_i$, only
    depends on the signs of $x_1-\zeta_{i,j}(y)$, $j=1,...,q_i$ ;
    \item  the zero set of $P_k$ coincide with the graphs of
    $\zeta_{i,j}$ ;
    \item the integers $m$, $q_i$, $deg(A_i)$, $deg(\zeta_{i,j})$ are bounded by
    a function of
    $\sum_k deg(P_k)$ and $d$.

\end{itemize}
\end{prop}

\begin{demo} \cite{Coste} Thm 2.3 p 112.
\end{demo}

From the above we deduce easily the following proposition :

\begin{prop}\label{pro3}
For all  semi-algebraic subsets $A\subset]0,1[^{d+1}$ , there exist
integers  $m$, $q_1,...,q_m$, a partition of $]0,1[^d$ into Nash
manifolds $A_1,...,A_m$ and Nash maps,
$\zeta_{i,1}<...<\zeta_{i,q_i}:A_i\rightarrow ]0,1[$, for all $1
\leq i \leq m$,  such that :
\begin{itemize}
    \item $A$ coincide with a union of slices of the two following forms
    $\{(x_1,y)\in ]0,1[\times A_i : \zeta_{i,k}(y)< x_1 <
    \zeta_{i,k+1}(y)\}$  and $\{(\zeta_{i,k}(y),y) :  y\in A_i\}$;
    \item the integers $m$, $q_i$, $deg(A_i)$, $deg(\zeta_{i,j})$ are bounded by a function of $deg(A)$ and $d$.
\end{itemize}
\end{prop}

For open semi-algebraic sets, we have the following result :

\begin{coro}\label{sc}
For all  semi-algebraic open subsets $A\subset]0,1[^{d+1}$, there exist integers $m$, $q_1,...,q_m$,
disjoint semi-algebraic open sets $A_1,...,A_m$ and Nash
maps, $\zeta_{i,1}<...<\zeta_{i,q_i}:A_i\rightarrow ]0,1[$, for all
$1 \leq i \leq m$,  such that :
\begin{itemize}
    \item $adh(A)$ coincide with a union of ``slices'' of the following form
$adh(\{(x_1,y)\in ]0,1[\times A_i : \zeta_{i,k}(y)< x_1
    <\zeta_{i,k+1}(y)\})$ ;
    \item the integers $m$, $q_i$, $deg(A_i)$, $deg(\zeta_{i,j})$ are bounded by a function of $deg(A)$ and $d$.
\end{itemize}
\end{coro}

In the following corollary, we reparametrize a semi-algebraic set with
Nash maps of bounded degree. 

\begin{coro}\label{1pas} (decomposition into cells) There exists $M:\mathbb{N}^3\rightarrow \mathbb{R}^+$, such that any semi-algebraic set  $A\subset ]0,1[^d$ admits a $(M)$-resolution.
\end{coro}

\begin{demo} We argue by induction on $d$. We note $P(d)$
the claim of the above corollary. $P(0)$ is trivial. Assume $P(d)$.

Let $A\subset ]0,1[^{d+1}$ be a semi-algebraic set  of maximum dimension
$l$. Proposition \ref{pro3} gives us integers  $m$, $q_1,...,q_m$, Nash manifolds $A_1,...,A_m\subset]0,1[^d$ and Nash maps,
$\zeta_{i,1}<...<\zeta_{i,q_i}:A_i\rightarrow ]0,1[$ such that :
\begin{itemize}
    \item $A$ coincides with an union of slices of the two following forms
    $\{(x_1,y)\in ]0,1[\times A_i : \zeta_{i,k}(y)< x_1 <
    \zeta_{i,k+1}(y)\}$  and $\{(\zeta_{i,k}(y),y) :  y\in A_i\}$;
    \item $m$, $q_i$, $deg(A_i)$ are bounded by a function of $deg(A)$ and $d$.
\end{itemize}

Let $1\leq i\leq m$. We apply the induction hypothesis to
$A_i\subset ]0,1[^d$ : there exist a resolution of $A_i$,
\emph{i.e.} an integer $N_i$ (bounded by a function of $deg(A_i)$ and  $d$,
therefore by a function of $deg(A)$ and $d$) and
Nash maps
$\phi_{i,1},...,\phi_{i,N_i}:]0,1[^l\rightarrow]0,1[^{d}$, such that
 $A_i=\bigcup_{p=1}^{N_i} \phi_{i,p}(]0,1[^l)$. Then we 
define
 $\psi_{i,k,p}:]0,1[^l\rightarrow ]0,1[^{d+1}$
as follows :
$\psi_{i,k,p}(x_1,y):=(x_1(\zeta_{i,k+1}(y)-\zeta_{i,k}(y))\circ\phi_{i,p}(x_1,...,x_{l})+
                                \zeta_{i,k}\circ\phi_{i,p}(x_1,...,x_{l}),\phi_{i,p}(x_1,...,x_{l}))$ for the first form and 
$\psi_{i,k,p}(x_1,y):=(\zeta_{i,k}\circ\phi_{i,p}(x_1,...,x_{l}),\phi_{i,p}(x_1,...,x_{l}))$ for the second form.\\

The $\psi_{i,k,p}$ are Nash triangular maps, such that
\begin{itemize}
\item the number of these parametrizations is bounded by $3\sum _{i=1}^m q_iN_i$
\item $A=\bigcup_{i,k,p}\psi_{i,k,p}(]0,1[^{l})$
\item $deg(\psi_{i,k,p})$ is bounded by a function of $deg(A)$ and $d$
(See Corollary 3)
\end{itemize}

Thus these  maps form a resolution of $A$.
\end{demo}

The following lemma is another application of the Tarski's principle  :

\begin{lem}\label{lem} Let $A\subset \mathbb{R}^d$ be a semi-algebraic
open set, $f:A\rightarrow \mathbb{R}^n$ a Nash map defined on $A$.
The partial derivatives of $f$ of all orders are also semi-algebraic maps
of degree bounded by a function of $deg(f)$, $d$ and $n$.
\end{lem}

\begin{demo} Apply corollary \ref{quant}. See \cite{BCR} p 29.
\end{demo}

\section{Proof of the Yomdin-Gromov Theorem}

First we show the following technical proposition, in which we work with
``smooth'' functions. Finally we explain how we reduce  the proof of the main theorem to this proposition.

\begin{prop}\label{pro} Let $A\subset ]0,1[^d$ be a semi-algebraic open set. Let $f_1,...,f_k:A\rightarrow
]0,1[$ be Nash maps   and let
$\alpha \in \mathbb{N}^d-\{0\}$.  There exists a sequence $(A_n)_{n\in
\mathbb{N}}\subset A^{\mathbb{N}}$ of semi-algebraic sets, such that\\

\begin{itemize}
  \item $a_n:=\sup_{x\in A} d(x,A_n) \xrightarrow[n\to + \infty]{} 0$  $(\#)$, where
$d(x,A_n)$ is the distance between $x$ and $A_n$ ;
  \item  $deg(A_n)$ is bounded by a function of $\max_i(deg(f_i))$, $|\alpha|$ and $d$ ;
  \item $(f_{i/A_n})_{i=1,...,k}$ admits a $(\C^{\alpha},1)$-resolution.
\end{itemize}

We will say that such a sequence  $(A_n)_{n\in\mathbb{N}}$ is
$\alpha$-adapted to $(f_i)_{i=1,...,k}$.
\end{prop}

The following corollary follows from the above proposition :

\begin{coro}\label{cor} There exists $M:\mathbb{N}^4\rightarrow \mathbb{R}^+$, such that for all integers $k\geq 1$, $d\geq 1$, multiindices $\alpha \in
 \mathbb{N}^d-\{0\}$,  any family  $(f_i:A\rightarrow
]0,1[)_{i=1,...,k}$ of Nash maps, where $A\subset
]0,1[^d$ is a semi-algebraic open set, admits
a $(\C^{\alpha},1,M)$-resolution.
\end{coro}

Now we show how Proposition 4,
 Corollary 6 and the Yomdin-Gromov theorem follow from the case $k=1$ of the proposition 4. In fact we show stronger results, which are used in the induction in the last section.

\begin{nota}\label{not} We consider the set $E$ of pairs
$(\alpha,d)$, where $d\in \mathbb{N}^{*}$ and $\alpha \in
\mathbb{N}^d-\{0\}$. The set $E$ is provided with the following order 
$\ll$ :\\

$(\beta,e)\ll(\alpha,d)$ iff $(e<d)$ or ($e=d$ and
$\beta\preceq\alpha$)\\

We will write :\\
 $P4(\alpha,d)$
the claim of the proposition 4 for all pairs $(\beta,e)$ with
$(\beta,e)\ll(\alpha,d)$.\\
$C6(\alpha,d)$ the claim of the corollary 6 for all pairs $(\beta,e)$ with
$(\beta,e)\ll(\alpha,d)$.\\
$YG(\alpha,d)$ the existence of a $(\C^{\beta},1)$ resolution for all Nash manifolds $A\subset [0,1]^{e}$, and for all pairs
$(\beta,e)\ll(\alpha,d)$.\\

\end{nota}

\begin{rem}
With the above notations, we have :  theorem 1 $\Longleftrightarrow YG(\alpha,d) \quad \forall (\alpha,d)\in E$.
\end{rem}

\begin{lem} \label{lemm}
The claim $C6(\alpha,d)$ [resp.  $P4(\alpha,d)$] for $k=1$
 implies the claim $C6(\alpha,d)$ [resp. $P4(\alpha,d)$] for all $k\in\mathbb{N}^{*}$.
 \end{lem}

\begin{demof}{Lemma 5 (Case of $C6(\alpha,d)$)} We argue by induction on $k$.

Assume that for $k$-families $g_1,...,g_k:B\rightarrow ]0,1[$
of Nash maps, with $B\subset ]0,1[^{d}$ a semi-algebraic open set, admit a
$(\C^{\alpha},1)$-resolution.

Let $f_1,...,f_{k+1}:A\rightarrow ]0,1[$ be
Nash maps, with $A\subset ]0,1[^{d}$ a semi-algebraic open set. According to the induction hypothesis, there exists
$(\phi_i)_{i=1,...,N}$ a $(\C^{\alpha},1)$-resolution of
$(f_1,...,f_k)$. According to $C6(\alpha,d)$ for $k=1$,  for each $i$ we can find
  $(\psi_{i,j})_{j=1,...,N_i}$ a $(\C^{\alpha},1)$-resolution  of 
$f_{k+1}\circ \phi_i$. According to Lemma \ref{lemcomp}, the maps
 $\phi_i\circ \psi_{i,j}$, of which the number is
 $\sum_{i=1}^N N_i$, are
 $(\C^{\alpha},K)$-Nash triangular maps,  as well as the maps $f_p\circ
\phi_i\circ\psi_{i,j}$ for all $1\leq p \leq k$ (with
$K=K(|\alpha|,d)$). Finally, for each $i$,
$(\psi_{i,j})_{j=1,...,N_i}$ being a $(\C^{\alpha},1)$-resolution
of $f_{k+1}\circ \phi_i$, the maps  $f_{k+1}\circ
\phi_i\circ\psi_{i,j}$  are
$(\C^{\alpha},1)$-Nash maps. Moreover, we have in a trivial way :
$adh(A)=\bigcup_{i,j}\phi_i\circ\psi_j([0,1]^d)$.
We conclude the proof for $C6(\alpha,d)$ thanks to Lemma \ref{lem3}.
\end{demof}

\begin{demof}{Lemma 5 (Case of $P4(\alpha,d)$)}

We adapt the above proof for $P4(\alpha,d)$ as follows. Let
$(A_n)_{n\in\mathbb{N}}$ be a sequence $\alpha$-adapted to
$(f_i)_{i=1,...,k}$. Hence, for all  $n\in \mathbb{N}$, there exists 
$(\phi_j^n)_{j=1,...,N_n}$ a $(\C^{\alpha},1)$ resolution
of $(f_{i/A_n})_{i=1,...,k}$.  For $n$, $j$, let 
$(A^{n,j}_p)_{p\in\mathbb{N}}$ be a sequence $\alpha$-adapted to
$f_{k+1}\circ\phi_j^n$. We use the following remark, which is an easy 
consequence of the compactness of $[0,1]^d$ :

\begin{rem}\label{adapte}
If $(A_n)_{n\in \mathbb{N}}$ is a sequence of subsets of $]0,1[^l$ satisfying 
$\sup_{x \in ]0,1[^l}d(x,A_n)\xrightarrow[n \to +\infty]{} 0$ and  $\phi:]0,1[^l\rightarrow ]0,1[^d$ is a continuous map \footnote{possibly not uniformly continuous}, then $\sup_{x \in \phi(]0,1[^l)}d(x,\phi(A_n))\xrightarrow[n \to +\infty]{} 0$.
\end{rem}  

According to the above remark, we can choose an integer $p_{j,n}$ for each 
$n\in \mathbb{N}$ and each $1\leq j \leq N_n$, such that 
$\sup_{x \in \phi^n_j(]0,1[^l)}d(x,\phi^n_j(A^{n,j}_{p_{j,n}}))<1/n$.
Now, let us show that $B_n:=\bigcup _{j=1}^{N_n}\phi_j^n(A^{n,j}_{p_{j,n}})$
defines a sequence $\alpha$-adapted to $(f_i)_{i=1,...,k+1}$.

 Observe that $B_n$ is a semi-algebraic set because each
 $\phi_j^n$ is a  semi-algebraic map and each $A_p^{n,j}$
are semi-algebraic sets. Moreover $N_n$,
$deg(\phi_j^n)$ and $deg(A^{n,j}_{p_{j,n}})$ and therefore $deg(B_n)$ are bounded by a function of  $\max_i(deg(f_i))$, $|\alpha|$ and $d$.
Finally, we check the ``density condition'' $(\#)$. 

$\sup_{x \in A}d(x,B_n)\leq \sup_{x \in A}d(x,A_n)+
\max_{j=1,...,N_n}(\sup_{x \in\phi^n_j(]0,1[^l)}d(x,\phi^n_j(A^{n,j}_{p_{j,n}}))
\leq a_n+1/n \xrightarrow[n \to + \infty]{} 0$. 
\end{demof}

\begin{demof}{Corollary 6 ($ P4((0,...,0,r+1),d) \Rightarrow C6((0,...,0,r),d)$.)}
According to Lemma \ref{lemm}, it is enough to consider a single Nash map 
$f:A\rightarrow
]0,1[$, where $A\subset
]0,1[^d$ is a semi-algebraic open set. According to $ P4((0,...,0,r+1),d)$,
 there exists a $(0,...,0,r+1)$-adapted sequence   $(A_n)_{n \in \mathbb{N}}$ 
to $f$. Let $(\phi_{i}^k)_{i\leq N_k}$ be a $(\C^{r+1},1)$-resolution of
$f_{/A_k}$. By hypothesis, $N_k$ is bounded by a function of $deg(A_k)$ and
 $r$ and thus by a function
of $deg(A)$ and $r$ ; consequently $(N_k)_{k\in \mathbb{N}}$ is a bounded sequence. By extracting a subsequence, we can assume
$N_k=N$, for all $k\in \mathbb{N}$. According to the Ascoli theorem, $B(r+1)^{2N}$ is a compact set in
$ B(r)^{2N}$, where $B(r)$ is the closed unit ball of $\mathcal{C}^r(]0,1[^d)$ (set of $\mathcal{C }^r$ maps on $]0,1[^d$ onto $\mathbb{R}$). By 
extracting again a subsequence, we can assume that for each  $i=1,...,N$
 $(\phi_{i}^n)_{n\in \mathbb{N}}$ converge
in $\|.\|_r$ norm to a
$(\C^r,1)$-Nash map, $\psi_i$. Obviously $f\circ\psi_i$ is a
$(\C^r,1)$-Nash map. One  only needs to see
$\bigcup _{i=1,...,N} \psi_i([0,1]^d)=adh(A)$. It is enough to show that
$A\subset \bigcup _{i=1,...,N} \psi_i([0,1]^d)$. We have
$\psi_i([0,1]^d)\subset adh(A)$, for all $i$, by convergence of $\phi_i^n$ to $\psi_i$. Let $x\in A$.
According to the ''density condition'' $(\#)$, there exists a sequence  $x_n\in
A_n$, such that $x_n\rightarrow x$. By extracting a subsequence, we can assume that there exist $1\leq i\leq N$ and a sequence
$(y_n\in [0,1]^d)_{n\in\mathbb{N}}$ such that $x_n=\phi_{i}^n(y_n)$.
  By the uniform convergence of $\phi_{i}^n$ to $\psi_i$, we have $\psi_i(y_n)
\rightarrow x$. We easily conclude that
$\bigcup _{i=1,...,N} \psi_i([0,1]^d)=adh(A)$.
Finally $(\psi_i)_{i\leq N}$ is a $(\C^r,1)$-resolution of $f$.
\end{demof}

\begin{demof}{Theorem 1 (  $C6(\alpha,d)\Rightarrow YG(\alpha,d+1)$)}
 
Under Proposition \ref{pro3}, it is enough to consider the two following special cases :

\begin{enumerate}
\item $A\subset ]0,1[^{d+1}$ is a semi-algebraic set of the form : $\{(x_1,y)\in ]0,1[\times A' : \eta(y)< x_1 <
    \zeta(y)\}$, where $A'\subset]0,1[^d$ is a semi-algebraic set of maximum dimension $e$ and $\eta, \zeta :
A'\rightarrow ]0,1[$  Nash maps, such that
$deg(\eta),deg(\zeta),deg(A')$ depend only on $deg(A)$ and $d$.
 By using a $\alpha$-resolution of $A'$  $(\phi_i:]0,1[^e\rightarrow ]0,1[^{d})_
{i=1,...,N}$ and by considering $\eta \circ \phi_i$ and $\zeta \circ \phi_i$, 
we can assume that $A'=]0,1[^{e}$, with $e\leq d$.  

Applying
 $C6(\alpha,d)$ to $(\zeta,\eta)$, there exists $(\phi_i)
_{i=1,...,N}$ a $(\C^{\alpha},1)$-resolution of $(\zeta,\eta)$.

For each  $i$, we define $\psi_i:]0,1[\times ]0,1[^{e}\rightarrow
]0,1[^{d+1}$ in the following way : $\psi_i(x,y)=(x(\zeta\circ\phi_i-\eta
\circ\phi_i)(y)+\eta\circ\phi_i(y), \phi_i(y))$. Then
$(\psi_i)_{i=1,...,N}$ is a $(\C^{\alpha},2)$-resolution of $A$.
We conclude the proof using Lemma \ref{lem3}. \\

\item $A$ is a semi-algebraic set of the form $\{(\zeta_{i,k}(y),y) :  y\in A'\}$. The decomposition into cells gives  (see Corollary \ref{1pas})  us  a resolution of $A$,  $(\phi_i:]0,1[^{l}\rightarrow ]0,1[^{d+1})_
{i=1,...,N}$, with $l<d+1$. We conclude the proof, by applying for each $i$,
$C6(\alpha,d)$ to the coordinates of $\phi_{i}$.
\end{enumerate}

\end{demof}

\section{Proof of Corollary 6 in dimension 1}

First we study the case of dimension $1$, where we can prove right away
Corollary 6. The case of dimension $1$ allows us to introduce simple ideas
of parametrizations, which will be adapted in higher dimensions.

The semi-algebraic sets of $]0,1[$ are the finite unions
of open intervals and points. So it's enough to prove the Corollary 6 for $A$ of the form  $]a,b[\subset]0,1[$. \\

\begin{demof}{$C6(1,1)$ (Case of the first derivative) }
Let $f:]a,b[\rightarrow ]0,1[$ be a $\C^0$-Nash map. \footnote{In dimension 1, a bounded Nash map (defined on a bounded intervall) is a $\C^0$-Nash map (See \cite{BCR} p 30)}
We cut the interval $]a,b[$ into a minimal number $N$ of subintervals
($J_k)_{k=1,...,N}$, such thatfor each $k$,  $\forall x\in J_k$,
$|f'(x)|\geq 1$ or $\forall x\in J_k$, $|f'(x)|\leq 1$.

The required bound  on $N$ results from Tarski's principle.\\
On each interval $J_k$, we consider the following parametrization $\phi$ of $adh(J_k)=[c,d]\subset [0,1]$ :
\begin{itemize}
    \item $\phi(t)=c+t(d-c)$ if $|f'|\leq 1$ and then we have $deg(\phi)=1$, $deg(f\circ\phi)=deg(f)$.
    \item $\phi(t)=f_{|[c,d]}^{-1}(f(c)+t(f(d)-f(c)))$ if $|f'|\geq 1$ and then we have $deg(\phi)=deg(f)$ (indeed
    $deg(f^{-1})=deg(f)$) and $deg(f\circ\phi)=1$.
\end{itemize}
\end{demof}

\begin{demof}{$C6(r,1)$ (Case of higher derivatives) } We argue by induction
on $r$ : assume $C6(r,1)$, with $r\geq 1$ and prove
$C6(r+1,1)$.\\

Let  $f:]a,b[\subset ]0,1[ \rightarrow ]0,1[$ be a $\C^0$-
Nash map. By considering $(f\circ\phi_i)_{i=1,...,N}$,
where $(\phi_i)_{i=1...N}$ is a $(\C^r,1)$-resolution of $f$ given by $C6(r,1)$,
we can assume that $f$ is a $(C^r,1)$-Nash map.\\

We divide the interval $]a,b[$ into a minimal number $n_i$ of subintervals
on which $|f^{(r+1)}|$ is either increasing or decreasing, \emph{ie}, the sign of $f^{(r+1)}f^{(r+2)}$ is  constant. Consider the case where  $|f^{(r+1)}|$ is decreasing, the increasing case being similar. We reparametrize those intervals from $[0,1]$ with linear increasing maps $\widetilde{\phi_i}$.
We define $f_i=f\circ\widetilde{\phi _i}$. Obviously $f_i$ is $C^r,1)$-Nash 
map and  $|f_i^{(r+1)}|$ is decreasing.  In the following computations, we note $f$ instead of $f_i$.\\

Setting $h(x)=x^2$, we have :
$$(f\circ h)^{(r+1)}(x)=(2x)^{r+1}f^{(r+1)}(x^2)+R(x,f(x),...f^{(r)}(x))$$
where $R$ is a polynomial depending only on $r$. Therefore
$$\forall x\in[0,1] \quad |(f\circ h)^{(r+1)}(x)|\leq |(2x)^{r+1}f^{(r+1)}(x^2)|+C(r),$$ where
$C(r)$ is a function of $r$.\\
Furthermore, we have \begin{equation} \label{equ1} x|f^{(r+1)}(x)|=\int_0^x
|f^{(r+1)}(x)|dt\leq | \int_0^x
f^{(r+1)}(t)dt|=|f^{(r)}(x)-f^{(r)}(0)|\leq2
\end{equation} thus \\
$$|(f\circ h)^{(r+1)}(x)|\leq C(r)+ 2\frac{(2x)^{r+1}}{x^2}\leq C(r)+2^{r+2}$$ Enfin
$deg(\widetilde{\phi_i}\circ h)=2$ and $deg(f\circ h)=2deg(f)$.
The claim concerning the integers $n_i$ results from the Tarski's principle. We conclude the proof of $C6(r+1,d)$ thanks to the lemma
3.

\end{demof}

\section{Proof of  Proposition 4}

Let us fix two integers $r\geq 2,c\geq 1$. In this section we show 
$P4((0,...,0,r-1),c)$
for $k=1$, as this implies the general case by  Lemma \ref{lemm}.\\

We argue by induction on the set $E_{rc}$ of pairs 
$(\alpha,d)$, where $d\in \mathbb{N}^{*},d\leq c$ and $\alpha \in
\mathbb{N}^d,
 |\alpha|\leq r+c-d$.  $E_{rc}$ is provided with the order $\ll$.\\

We assume now that $P4(\alpha,d)$ is checked and we distinguish three cases
depending on the values of the pair $(\alpha,d)$ :\\
\\
\begin{tabular}{|c|}
  \hline
  \emph{Increase of the dimension : $P4((0,...,0,r+c-d),d)\Rightarrow P4((1,0,...,0),d+1)$}  \\
  \hline
\end{tabular}

\begin{demo}

\begin{claim} \label{r2} It is enough to show the result for Nash maps $f:]0,1[^{d+1}\rightarrow ]0,1[$.
\end{claim}

\begin{demof}{Claim \ref{r2}}

Let $f:A\subset ]0,1[^{d+1} \rightarrow ]0,1[$ a 
Nash map,  defined on a semi-algebraic open set of
$\mathbb{R}^{d+1}$.

Consider a resolution
$(\phi_i:[0,1]^{d+1}\rightarrow [0,1]^{d+1})_{i=1,...,N}$ of $A$ given by Lemma \ref{1pas}. If $(A_n)_{n\in \mathbb{N}}$ is an adapted sequence to $( f\circ\phi_i,\phi_i)$ and $(\psi^{i,n}_j)_{j=1,...,N_{i,n}}$ a $C^{(1,0,...,0)}$ resolution of $( f\circ\phi_{i /A_n},\phi_{i /A_n})$, then under Remark \ref{adapte}, the sequence 
$(B_n)_{n\in \mathbb{N}}$, defined as follows  
$B_n=\bigcup_{i=1,...,N} \phi_i(A_n)$, is an adapted sequence to $f$ with 
$(\phi_i\circ\psi^{i,n}_j)_{i,j}$ as a resolution of $f_{/B_n}$.      
\end{demof}

We work on $A_n=]1/n,1-1/n[^{d+1}$ in order to ensure that $f$ extends continuously on $adh(A_n)$. For simplicity, we note $A$ instead of $A_n$.

We consider the following semi-algebraic open sets :
$A_+=int(\{x\in A , |\partial_{x_1} f(x)|>1\})$ and $A_-=int(\{x\in
A,|\partial_{x_1} f(x)|\leq1\})$. We have
 $adh(A)=adh(A_+)\bigcup adh(A_-)$. Obviously $adh(A_+)\bigcup adh(A_-)\subset adh(A)$. Let show
 $A \subset adh(A_+)\bigcup adh(A_-)$. Let $x\in A$. If $d(x,A_n^+)=0$, then $x\in adh(A_+)$ ; if not, as $A$ is an open set, there exists
 $r>0$, such that the ball $B(x,r)\subset A\bigcap A_+^c\subset \{x\in A,|\partial_{x_1} f(x)|\leq1\}$ and thus
 $x\in A_-$.    \\

According to $P4((0,...,0,r+c-d),d)\Rightarrow P4((0,...,0,2),d)
\Rightarrow C6((0,...,0,1,d))\Rightarrow YG((0,...,0,1),d+1)$, there exist
 $(\C^{(1,0,...0)},1)$-Nash triangular maps
$(\phi_{j})_{1\leq j\leq N}$ such that $adh(A_-)=\bigcup_{1\leq j\leq
N_-}\phi_j([0,1]^d)$ and such that $N_-$, $deg(\phi_{j})$ are
bounded by a function of $deg(A_-)$, and thus by a function of $deg(f)$ (according to the lemma 4 and the corollary \ref{cor}). We have $|\partial_{x_1}(f\circ
\phi_{j})|\leq 1$, so the maps
$\phi_i$ can be used to build a resolution of $f$.\\

For $A_+$, we consider the inverse of $f$. Observe first, that according to the
corollary \ref{sc}, we can assume that $A_+$ is a slice of the following form $\{(x_1,y)\in ]0,1[ \times
A_+' \ : \zeta(y)< x_1< \eta(y)\}$, where $A_+'\subset ]0,1[^d$ is
a semi-algebraic open set of $\mathbb{R}^d$ and
$\zeta,\eta:A_+'\rightarrow ]0,1[$ are Nash maps.

Define $D_+=\{(f(x_1,y),y) : (x_1,y)\in A_+\}$. We define 
$g:A_+\rightarrow D_+$,  $g(x_1, y_1)=(f(x_1,y),y))$. This map $g$
 is a local diffeomorphism, by the local inversion theorem.
 Moreover, $g$ is one to one, because $g(x_1,y)=g(x'_1,y')$ implies
 $y=y'$, and $f(x_1,y)=f(x'_1,y)$ implies $x_1=x'_1$, because
$|\partial_{x_1} f(x)|\geq 1$ for $x\in A_+$.  The map $g$ extends to
$g:adh(A_+)\rightarrow adh(D_+)$, a homeomorphism, since
 $f$ is continuous on $adh(A)$ (Recall that we note $A:=A_n$).

Observe that $D_+$ is a semi-algebraic open set of $\mathbb{R}^{d+1}$.
On $D_+$ we define $\phi$ :
$\phi(t,u):=g^{-1}(t,u)=(f(.,u)^{-1}(t),u)$. The Nash map $\phi:D_+ \rightarrow
A_+$
is  triangular  and $deg(\phi)=deg(f)$.
Define $\phi (t,u)=(x_1,y)$. We compute  : $$ D\phi(t,u)=\left(%
\begin{array}{cc}
                         \frac{1}{\partial_{x_1}f(x_1,y)} &
                         -\frac{1}{\partial_{x_1}f}\nabla_{y} f(x_1,y) \\
                         0 & Id \\
                       \end{array}%
                       \right)$$

As $(x_1,y)\in A_+$, we have $|\partial_{x_1}\phi|\leq 1$. Furthermore, 
we check
$$f\circ\phi(t,u)=t.$$

Therefore, $\phi$ and $f\circ\phi$ are
$(\C^{(1,0,...,0)},1)$-Nash triangular maps. In order to obtain a resolution, we apply again $YG((0,0,...,0,1),d+1)$ to $adh(D_+)$. That gives a $(C^{(1,0,...,0)},1)$
-Nash triangular parametrization $\psi_{j}:[0,1]^{d+1}\rightarrow adh(D_+)$,
$j\leq N_+$, such that $N_+$,
 $deg(\psi_{j})$ are bounded by a function of $deg(D_+)$, thus by a function of $deg(f)$. Moreover 
$$|\partial_{x_1}(\phi\circ\psi_j)|=|\partial_{x_1}(\phi)|.|\partial_{x_1}(\psi_j^1)|\leq 1 $$ because $\psi_j$ is triangular and
$$|\partial_{x_1}(f\circ \phi \circ \psi_{j})|=|\partial _{x_1}
\psi_{j}^1|\leq 1,$$ where $\psi_{j}:=(\psi^1_j,...,\psi^{d+1}_j)$.
The following parametrizations $\phi\circ\psi_j:[0,1]^{d+1} \mapsto [0,1]^{d+1}$
are therefore $(\C^{(1,0,...,0)},1)$-Nash triangular maps
such that :
\begin{itemize}
    \item $adh(A_+)=\bigcup_{j=1}^{N_+}
    \phi\circ\psi_j([0,1]^{d+1})$ ;
    \item each $f\circ\phi\circ\psi_j$ is a
    $(\C^{(1,0,...,0)},1)$-Nash map ;
    \item  $deg(\phi\circ\psi_j)$, $deg(f\circ\phi\circ\psi_j)$ are bounded by a function of $|\alpha|,d$,
    and $deg(f)$ (See Corollary 3).
\end{itemize}

Finally, we combine the maps $\phi_1,...,\phi_{N_-}$ with the maps $\phi\circ\psi_1,...,\phi\circ\psi_{N_+}$, so that we obtain a  $(\C^{(1,0,...,0)},1)$-resolution of $f$. The bound on the number of parametrizations is the result of the bounds
 on  $N_-$ and
$N_+$ from the Yomdin-Gromov theorem and of the bounds from the proposition \ref{pro1}.

\end{demo}

\begin{tabular}{|c|}
  \hline
  \emph{ Increase of the derivation order : $P4((0,...,0,s),d)\Rightarrow P4((s+1,0,...,0),d)$} pour
$s<r+c-d$ \\
  \hline
\end{tabular}\\
Until the end, $C(|\alpha|,d)$ are functions of $|\alpha|$ and $d$.

\begin{demo}

In this case, we adapt the proof in dimension $1$. We begin with a remark
 similar to the previous Claim \ref{r2}. 

\begin{claim} \label{rr}
It is enough to show the result for 
$(\C^s,1)$-Nash maps $f:A=]0,1[^d\rightarrow ]0,1[$.
\end{claim}

\begin{demof}{claim \ref{rr}}
Let $f:A\subset ]0,1[^{d+1} \rightarrow ]0,1[$ a 
Nash map,  defined on a semi-algebraic open set of
$\mathbb{R}^{d+1}$.
By applying $P4((0,...,0,s),d)$  to $f$, we obtain a $(\C^s,1)$-resolution
$(\phi_i^n)_{i=1,...,N_n}$ of $f_{/A_n}$, with $A_n$  an adapted sequence.
 We conclude by applying  $P4((s+1,0,...,0),d)$ to the family of
 $(\C^s,1)$-Nash maps $(f\circ\phi_i^n,\phi_i^n)$, and by applying remark \ref{adapte}.
\end{demof}

Let $f:A=]0,1[^d\rightarrow ]0,1[$ be a
$(\C^s,1)$-Nash map. 

 We cut up $]0,1[^d$ according to the sign  $\frac{\partial ^{s+1} f}{\partial
x_1^{s+1}}\frac{\partial ^{s+2} f}{\partial x_1^{s+2}}$, and we 
assume (See corollary \ref{sc}) that $A$ is a slice of the following
form $\{(x_1,y) \in ]0,1[\times A' \ : \zeta(y)<x_1<\eta(y)\}$,
where $A'\subset]0,1[^{d-1}$ is a semi-algebraic open set and 
$\zeta, \eta :A'\rightarrow ]0,1[$ are Nash maps.

Applying the estimate (\ref{equ1}) obtained in section 6
to the  function $x_1\mapsto\frac{\partial ^{s+1}
f}{\partial x_1^{s+1}}(x_1,y)$ (we fix $y$), we get

\begin{equation} \label{equ2} |\frac{\partial ^{s+1} f}{\partial
x_1^{s+1}}(x_1,y)|\leq \frac{2}{|x_1-\zeta(y)|}
\end{equation} or
\begin{equation}\label{equ3}|\frac{\partial ^{s+1} f}{\partial
x_1^{s+1}}(x_1,y)|\leq \frac{2}{|x_1-\eta(y)}|,
\end{equation} according to the sign of $\frac{\partial ^{s+1}
f}{\partial
x_1^{s+1}}\frac{\partial ^{s+2} f}{\partial x_1^{s+2}}$. \\

The induction hypothesis $P4((0,...,0,s),d)$ implies 
$P4((0,...,0,s+2),d-1)$, because $(0...0,s+2),d-1)\ll
 ((0,...,0,s),d))$ and  $P4((0,...,0,s+2),d-1)$ implies
$C6((0,...,0,s+1),d-1)$. Apply $C6((0,...,0,s+1),d-1)$ to
$(\zeta,\eta)$ : there exist $(\C^{s+1},d-1)$-Nash
triangular maps $h:[0,1]^{d-1}\rightarrow [0,1]^{d-1}$, of which the images 
cover $adh(A')$, such that $\zeta \circ h$ and $\eta \circ h$ are 
$(\C^{s+1},d-1)$-Nash maps. Define
$\psi:[0,1]\times[0,1]^{d-1}\rightarrow adh(A)$,
$$\psi(v_1,w)=(\zeta\circ h(w).(1-v_1^2)+\eta\circ h(w).v_1^2,h(w))$$
The map $\psi$ is triangular and
$\|\psi\|_{s+1}\leq 2$.\\

In the new coordinates $(v_1,v_2,...,v_{d})$, the previous estimates
 (\ref{equ2}) and (\ref{equ3})
become, with $w=(v_2,...,v_{d})$ : $$|\frac{\partial ^{s+1}
f}{\partial x_1^{s+1}}(\psi(v_1,w))|\leq \frac{2}{v_1^2|\eta\circ
h(w)-\zeta\circ h(w)|}$$

Moreover, $\frac{\partial^{s+1}(f\circ\psi)}{\partial
v_1^{s+1}}(v_1,w)=(2v_1)^{s+1}(\eta\circ h(w)-\zeta\circ
h(w))\frac{\partial ^{s+1} f}{\partial
x_1^{s+1}}(\psi(v_1,w))+R(\eta\circ h(w)-\zeta\circ
h(w),v_1,(\frac{\partial ^{k} f}{\partial
x_1^{k}}(\psi(v_1,w)))_{k\leq s})$, where $R$ is a polynomial, which depends
only on $s$. The first part is less than $2^{s-1}$. Consider the second part. 
The map  $f$ is a $(\C^{s},1)$-Nash map, therefore $|\frac{\partial ^{k} f}{\partial x_1^{k}}|\leq
1$, for $k\leq s$ ; thus $|R(\eta\circ h(w)-\zeta\circ
h(w),v_1,(\frac{\partial ^{k} f}{\partial
x_1^{k}}(\psi(v_1,w)))_{k\leq s})|$ is bounded by a function of $s$, and
 therefore $|\frac{\partial^{s+1}(f\circ\psi)} {\partial v_1^{s+1}}|\leq
C(s,d)$. According to lemma \ref{comp}, the derivatives of lower order than $s$ of $f\circ\psi$
are also bounded by a function of $s$. Using Lemma \ref{lem3},
 we can assume that $\psi$ is a $(\C^{s+1},1)$-Nash map and $f\circ\psi$ is a 
$(\C^{(s+1,0,...,0)},1)$-Nash map. \\

 By Lemma 4, Proposition 2 and
$C6((0,...,0,s+1),d-1)$ the number of
parametrizations $h$ and their degree are also bounded by such a function.
It follows that the total number of parametrizations $\psi$ is bounded by a function of $d$ and of $deg(f)$. We conclude using Corollary \ref{cor3}, that the same holds for the degree of the parametrizations $\psi$.

\end{demo}

\begin{tabular}{|c|}
  \hline

  \emph{Control of the following derivative : $P4(\alpha,d)\Rightarrow P4(\alpha+1,d)$} with
$\alpha\neq (0,...,0,s+1)$\\
  \hline
\end{tabular}

\begin{demo}
According to the Claim \ref{rr}, we can assume that $f:]0,1[^d\rightarrow ]0,1[$ is a $(\C^{\alpha},1)$-Nash map.

Define  $A_n=]1/n,1-1/n[^{d-1}$ and $b_n=1-2/n$.
According to the Tarski's principle
(See Corollary \ref{quant}),
$B=\{(x_{1},y)\in adh(A_n) \ : \ |\frac{\partial
^{\alpha +1} f}{\partial x^{\alpha
+1}}(x_{1},y)|=sup_{t\in[1/n,1-1/n]}(|\frac{\partial ^{\alpha +1}
f}{\partial x^{\alpha +1}}(t,y)|)\}$ is a semi-algebraic set of degree bounded
  by a function of $deg(f)$ and $s$. .
We have introduced the concept of adapted sequence, so that
the $sup$ above is bounded (recall that $f$ is not supposed
analytic in a neighbourhood of $A$). According to Proposition \ref{pro3},
$B$ is covered by sets $(B_i)_{i=1,...,N}$, $B_i=\{
(x_1,y)\in]0,1[\times B_i' : \gamma_i(y)< x_1<
\Delta_i(y)\}$ or $B_i=\{
(\sigma(y),y)\in B_i'\}$ , where $B'_i\subset ]1/n,1-1/n[^{d-1}$ are semi-algebraic 
sets of $\mathbb{R}^{d-1}$, such that $\bigcup_{i=1}^N
B'_i=]1/n,1-1/n[^{d-1}$ and where
$\sigma_i,\gamma_i,\Delta_i:B_i'\rightarrow ]0,1[$ are
Nash maps. In the first case, we set $\sigma_i:=1/2(\Delta_i +\gamma_i)$. 
Afterwards, we consider only the open sets $B'_i$. Observe that for these sets we have 
$\bigcup adh(B'_i)=[1/n,1-1/n]^{d-1}$.

We check thanks to the Tarski's principle and 
the proposition \ref{pro} that $N$ and the degree of
$\sigma_i$ are bounded by a function of $deg(f)$ and $s$.
Define $g_i(y)=\frac{\partial ^{{(\alpha +1)}_1} f}{\partial
x_1^{{(\alpha +1)}_1}}(\sigma_i(y),y)$ with $y\in adh(B'_i)$,
where $(\alpha+1)_{i}$ represent the $i^{th}$ coordinate of $\alpha+1$.
The induction hypothesis $P4(\alpha,d)$ implies
$P4((0,...,0,|\alpha|+1),d-1)$ and thus $C6((0,...,0,|\alpha|),d-1)$, which 
applied to $\sigma_i$ et $g_i$ gives
$(\C^{|\alpha|},1)$-Nash triangular maps
$h_{i,k}:[0,1]^{d-1}\rightarrow [0,1]^{d-1}$,  such that $g_i\circ
h_{i,k}$ and $\sigma_i\circ h_{i,k}$ are $(\C^{|\alpha|},1)$-Nash and such that
 $\bigcup_k h_{i,k}([0,1]^{d-1})=adh(B_i')$.\\

Then,
$$\frac{\partial ^{((\alpha +1)_2,...,(\alpha+1)_d)}(g_i\circ h_{i,k})}{\partial x^{((\alpha +1)_2,...,(\alpha +1)_d)}}(y)
=\frac{\partial ^{\alpha +1} f}{\partial x^{\alpha
+1}}(\sigma_i\circ h_{i,k}(y),h_{i,k}(y))\times (\frac{\partial
h_{i,k}}{\partial_{x_2}})^{{(\alpha +1)}_2}...(\frac{\partial
h_{i,k}}{\partial x_{d}})^{{(\alpha +1)}_{d}}+R$$

where $R$ is a polynomial of  derivatives of order $\preceq
\alpha$, and of the derivatives of $h_{i,k}$ and $\sigma_i\circ
h_{i,k}$ of order less than $|\alpha|$, $R$  depending only on $\alpha$.
The map $h_{i,k}$ is a $(\C^{|\alpha|},1)$-Nash map and  by hypothesis $f$ 
is a $(\C^{\alpha},1)$-Nash map, so that we have $|R|<C(|\alpha|,d)$.\\
After all $g_i\circ h_{i,k}$ is a $(\C^{|\alpha|},1)$-Nash map. Hence we have \\
$$|\frac{\partial ^{\alpha+1} f}{\partial x^{\alpha +1}}(\sigma_i\circ
h_{i,k}(y),h_{i,k}(y)) (\frac{\partial
h_{i,k}}{\partial_{x_2}})^{{(\alpha +1)}_2}...(\frac{\partial
h_{i,k}}{\partial x_{d}})^{{(\alpha +1)}_{d}}|\leq|\frac{\partial
^{((\alpha +1)_2,...,(\alpha +1)_d)}(g_i\circ h_{i,k})}
{\partial x^{((\alpha +1)_2,...,(\alpha +1)_d}}|+|R|<C(|\alpha|,d)$$\\

Define $\phi_{i,k}:[0,1]^{d}\rightarrow [0,1]^d$ by :
$$\phi_{i,k}(x_1,y)=(1/n + b_n x_1,h_{i,k}(y))$$ 

$\phi_{i,k}$ is a $(\C^{\alpha+1},1)$-Nash triangular map. We
check the two following points :
\begin{itemize}
    \item $\frac{\partial ^{\alpha +1} (f\circ\phi_{i,k})}{\partial
    x^{\alpha +1}}=\frac{\partial ^{\alpha +1}
f}{\partial x^{\alpha +1}}(1/n+b_nx_1,h_{i,k}(y))\times
(b_n)^{(\alpha+1)_1}
 (\frac{\partial
h_{i,k}}{\partial_{x_2}})^{{(\alpha +1)}_2}...(\frac{\partial
h_{i,k}}{\partial x_{d}})^{{(\alpha +1)}_{d}}+S$, where $S$ is a polynomial
 of the derivatives of $f$ of order $\beta\preceq \alpha$ (because $h_{i,k}$
is triangular) and of the derivatives of $h_{i,k}$ of order less than
 $|\alpha|$, $S$ depending only on $\alpha$. From above we deduce that
$|S|<C(|\alpha|,d)$. \\

Moreover by definition of $\sigma_i$, $|\frac{\partial ^{\alpha
+1}f}{\partial x^{\alpha +1}}(1/n+b_n x_1,h_{i,k}(y))\times
(\frac{\partial h_{i,k}}{\partial_{x_2}})^{{(\alpha
+1)}_2}...(\frac{\partial h_{i,k}}{\partial x_{d}})^{{(\alpha
+1)}_{d}}|\leq |\frac{\partial ^{\alpha +1}f}{\partial x^{\alpha
+1}}(\sigma_i \circ h_{i,k}(y),h_{i,k}(y))\times (\frac{\partial
h_{i,k}}{\partial_{x_2}})^{{(\alpha +1)}_2}...(\frac{\partial
h_{i,k}}{\partial x_{d}})^{{(\alpha +1)}_{d}}|<C(|\alpha|,d)$,\\
thus $|\frac{\partial ^{\alpha +1} (f\circ\phi_{i,k})}{\partial
    x^{\alpha +1}}|\leq|\frac{\partial ^{\alpha +1}
f}{\partial x^{\alpha +1}}(1/n+ b_n x_1,h_{i,k}(y))\times
(\frac{\partial h_{i,k}}{\partial_{x_2}})^{{(\alpha
+1)}_2}...(\frac{\partial
h_{i,k}}{\partial x_{d}})^{{(\alpha +1)}_{d}}|+|S|<C(|\alpha|,d)$\\
    \item finally for $\beta\preceq\alpha$, in the expression $\frac{\partial ^\beta
(f\circ\phi_{i,k})}{\partial
    x^\beta}$ take part only the derivatives of  $f$ of order
    $\preceq\alpha$, again because of the triangularity of $h_{i,k}$.
    Hence $|\frac{\partial ^\beta (f\circ\phi_{i,k})}{\partial
    x^\beta}|<C(|\alpha|,d)$.\\
\end{itemize}

The lemma \ref{lem3} gives us a $(\C^{\alpha},1)$-resolution of
$f_{/A_n}$. 

\end{demo}

\medbreak
This work is part of the author's Master thesis (Master 2 at Universite Paris-Sud) with the supervision of J. Buzzi.

David Burguet \\ CMLS - CNRS UMR 7640\\
Ecole polytechnique \\ 91128 Palaiseau Cedex \\ France\\
burguet@math.polytechnique.fr


\newpage



\begin{thebibliography}{99999}

\bibitem{BR}
R. Benedetti, J.-L. Riesler, {\it Real algebraic geometry and semi-algebraic
sets}, Hermann (1990).
\bibitem{BD}
M. Boyle, T. Downarowicz, {\it The entropy theory of symbolic
extension}, Invent. Math. 156 (2004), no.1, 119-161 .
\bibitem{BFF}
M. Boyle, D. Fiebig, U. Fiebig, {\it Residual entropy, conditional
entropy and subshift covers}, Forum Math. 14 (2002), 713-757.
\bibitem{Bu0}
J. Buzzi, {\it Ergodic and topological complexity of dynamical
systems}, Course given during the Research trimester Dynamical Systems, 2002, Pise.
\bibitem{Bu}
J. Buzzi, {\it Intrinsic ergodicity of smooth interval maps}, Israel
J.Math.100 (1997),  125-161.
\bibitem{BCR}
J. Bochnak, M. Coste, M.F. Roy, {\it Real semi-algebraic geometry},
Ergebnisse der Mathematik und ihrer Grenzgebiete (3), 36,
Springer-Verlag (1978).
\bibitem{C&Y}
G. Comte, Y. Yomdin, {\it Tame Geometry with applications in smooth
Analysis}, L.N.M. 1834, Springer (2004).
\bibitem{Coste}
M. Coste, {\it Ensembles semi-alg\'ebriques}  109-138, in : L.N.M. 959
{\it Géométrie algébrique réelles et formes quadratiques}, Springer
(1982).
\bibitem{Cow}
W.Cowieson, L.-S. Young, {\it SRB measures as Zero-noise limits}, Ergodic Theory and dynamical
systems 25 (2005),  1115-1138.
\bibitem{New1}
T. Downarowicz, S. Newhouse, {\it Symbolic extensions and smooth dynamical systems} Invent. Math. (on line) (2004).

\bibitem{Gr}
M. Gromov, {\it Entropy, homology and semi-algebraic geometry},
S\'eminaire Bourbaki 663 (1986).
\bibitem{New}
S. Newhouse, {\it Continuity properties of the entropy}, Annals of
Math. 129 (1989),  215-237.
\bibitem{w}
A. Wilkie, http://www.newton.cam.ac.uk/webseminars/pg+ws/2005/maa/0203/wilkie.
\bibitem{Yoma}
Y. Yomdin, {\it Volume growth and entropy}, Israel J.Math. 57 (1987), 285-300.
\bibitem{Yomb}
Y. Yomdin, {\it $\mathcal{C}^r$-resolution}, Israel J.Math. 57 (1987), 301-317.


\end{thebibliography}
\end{document}